\documentstyle[11pt, a4, amstex, amscd, %
epsfig, euscript,amssymb, twoside]{article}

\newtheorem{theorem}{Theorem}[section]
\newtheorem{lemma}[theorem]{Lemma}

\numberwithin{equation}{section}

\def\d{\delta}
\def\D{\Delta}
\def\e{\varepsilon}

\def\a{\alpha}
\def\b{\beta}
\def\uu{{\frak{u}}_q}
\def\oq{\overline{\Q}}
\def\t{\underline{t}}
\def\fg{{\frak{g}}}
\def\ol{\overline}
\def\cc{{\cal{C}}}
\def\u{{{\frak{u}}_q{(\cc)}}}
\def\s{{\sum\limits_{x\in\z_n^t}}}

\def\qqbin#1{
\left(\begin{array}{c}
    #1\\
    u\\
    \end{array}\right )
}

\def\qbin#1#2{
\left(\begin{array}{c}
    #1\\
    #2\\
    \end{array}\right )
}
\def\bin#1#2{
\bigg[\begin{array}{c}
    #1\\
    #2\\
    \end{array}\bigg]
}

\def\tsum#1{\sum\limits_{#1}}

\def\otm{\otimes}

\def\z{{\mathbb{Z}}}
\def\Q{{\mathcal{Q}}}

\newenvironment{proof}{{\it\noindent Proof.}}{\hfill $ \square $ \vskip 4mm}

\date{}

\begin{document}
\pagestyle{myheadings} \markboth{\centerline{\sc h.-l. huang, s.
yang}}{ \centerline{\sc quantum groups and double quiver
algebras}}

\title{\sc quantum groups
and double quiver algebras}
\author{{Hua-Lin Huang}\\
\small\sl Mathematics Section,
The Abdus Salam International Centre for Theoretical Physics
\\
\small\sl Strada Costiera 11, 34100 Trieste, Italy
\\ \footnotesize\sl and \\
\small\sl Department of Mathematics, University of Science and
Technology of China, \\
\small \sl Hefei 230026, Anhui, P. R.
China \\
\small\it hualin@@ustc.edu.cn
\\
\\
\and
 {Shilin Yang\footnote{The  author is
partially supported by the National Science Foundation of China (Grant. 10271014) and
Natural Science Foundation of Beijing City(Grant. 1042001)\newline
\ \qquad {\rm 2000} {\it Mathematics Subject Classification}.
81R50, 16W30. }} \\
 \small\sl College of Applied Sciences,
Beijing University of Technology\\
\small\sl
Beijing 100022, P. R.  China \\
\small\it slyang@@bjut.edu.cn} \maketitle

\begin{abstract}
For a finite dimensional semisimple Lie algebra $\fg$
and a root $q$ of unity in a field $k,$  we associate to these
data a double quiver $\ol{\Q}.$ It is shown that a restricted
version of the quantized enveloping algebras $U_q(\fg)$ is a
quotient of the double quiver algebra $k\ol{\Q}.$
\end{abstract}

\section*{Introduction}
Let $U_q(\fg)$  be the Drinfeld-Jimbo quantum group,
which is a deformation of
the universal enveloping algebra of a finite dimensional
semisimple Lie algebra $\fg.$  In the generic case, i.e. the
parameter $q$ is not a root of unity, several models have been
raised to realize it. For example, the Ringel-Hall algebra
approach is one successful model among them, see \cite{r, g, x}.
The case where $q$ is a root of unity is of
an particular interest since it is related with the modular
representation theory. It is remarkable that a finite dimensional
Hopf algebra, so-called restricted version of $U_q(\fg)$ arose
naturally when Lusztig considered this class of quantum group at roots of
unity. See \cite{l1,l2,l3}.  Just as the generic case, it is interesting
 to find an algebra method to realize $U_q(\fg),$
or its restricted version. In this sense, Cibils \cite{c1, c2}
found that some quotient of a particular path algebra
can realize the positive part of the
restricted quantized enveloping algebra corresponding to $\fg$.
Gordon in \cite{ig} extended several results in \cite{c2}.
We remark that  the restricted form of
$U_q({{\frak{sl}}_2})$   is realized  in the paper \cite{yang} by applying the deformation of
preprojective algebra introduced  in \cite{cbh}.

Let $k$ be a field, $\ell=n$ if $n$ is odd and $\ell=n/2$ if $n$
is even, where $n\ge 5$.
Let $\fg$ be a finite dimensional semisimple simple laced Lie algebra
and $q$ an $n$-th primitive root  of unity.  By $\cc$ we denote
the associated Cartan matrix of $\fg.$  For the Cartan matrix
$\cc=(a_{ij})_{t\times
t}$, there is an associated
quantum algebra $\u$, which by definition, is an associative
$k$-algebra with generators $K_i, \ E_i, \ F_i$ for $1 \le i \le
t$, subjecting to the relations
\begin{eqnarray}
     &&K_i^n=1, \ K_iK_j=K_jK_i; \label{eqn3-1}\\
    &&K_iE_jK_i^{-1}=q^{a_{ij}}E_j, \ K_iF_jK_i^{-1}=q^{-a_{ij}}F_j; \label{eqn3-2}\\
    && E_iF_j-F_jE_i=\delta_{ij}\frac{K_i-K_i^{-1}}{q-q^{-1}};\label{eqn3-3}\\
    &&E_i^\ell=0,\ F_i^\ell=0; \label{eqn3-4}\\
    &&\sum_{s=0}^{1-a_{ij}} (-1)^s \bin{1-a_{ij}}{s}_q
    E_i^{1-a_{ij}-s} E_j E_i^s=0, \textrm{ if } i\ne j; \label{eqn3-5}\\
    &&\sum_{s=0}^{1-a_{ij}} (-1)^s
    \bin{1-a_{ij}}{s}_q F_i^{1-a_{ij}-s} F_j F_i^s
    =0, \textrm{ if } i\ne j, \label{eqn3-6}
\end{eqnarray}
where $[m]_x=\frac{x^m-x^{-m}}{x-x^{-1}}$, $[m]!_x=[m]_x\cdots [1]_x,$
$[0]_x!=1$,  and
$\bigg[\begin{array}{c}
    {\small m}\\
    \small s\\
    \end{array}\bigg]_x=\frac{[m]!_x}{[s]!_x [m-s]!_x}$ for an indeterminate $x$.
The algebra $\u$  is a Hopf algebra, where the comultiplication, counit, and
antipode are given as follows
\begin{eqnarray*}
\D(K_i) &=& K_i \otimes K_i, \\ \D(E_i) &=& E_i \otimes 1 + K_i
\otimes E_i,\\ \D(F_i) &=& F_i \otimes K_i^{-1} + 1 \otimes F_i,
\\ \e(K_i)&=&1, \ \e(E_i)=\e(F_i)=0,
\\ S(K_i)&=& K_i^{-1}, \ S(E_i)=-K_i^{-1} E_i, \ S(F_i)=-F_i K_i.
\end{eqnarray*}
It would be noted that the restricted quantum group in the sense of Luzstig
is some quotient of $\u$ \cite{l2, ig}.

 As a continuation of \cite{yang},  the aim of this paper is to
construct the quantum group $\u$ from a quiver up to Hopf algebra
isomorphism for an arbitrary finite dimensional semisimple simply
laced Lie algebra $\fg$. Accordingly, we first define a double
quiver $\ol{\Q}$ associated to $\cc$ and $q$,  and we get the double
quiver algebra $k\ol{\Q}.$  Then we
describe an algebra $\Pi^\cc$ obtained from the  algebra
$k\ol{\Q},$ which is actually a Hopf algebra. Finally, we
construct  a quotient $\uu^\cc$ of $\Pi^\cc$, which inherits the
Hopf algebra structure. We will show that $\u$ is isomorphic to
$\uu^\cc$
as Hopf algebras.
It is mentioned that the results in \cite{c2} are extended.
If $\fg$ is the three dimensional simple
Lie algebra, $\Pi^\cc$ is just a deformation of preprojective
algebra in \cite{cbh} and $\uu^\cc$ is the algebra ${\frak u}_q$
described in \cite{yang}.

\section{Preliminaries}

In this section, we list some notations for the convenience of the
statement.

Let $k$ be a fixed field. Let us fix an integer $n \ge 5$, an
$n$-th primitive root $q$ of unity and a positive definite symmetric Cartan
matrix $\cc=(a_{ij})_{t\times t}.$  By $c^i$ we denote the $i$-th
column of $\cc$, by $\z_s$  the cyclic group $\frac{\Bbb Z}{s \Bbb
Z}$ for any positive integer $s$, and by $\t$  the set $\{1, 2,
\cdots, t\}$.

For $\a=(a_1, \cdots, a_t), \b=(b_1, \cdots, b_t)\in \z_n^t$, we
denote $\a\cdot \b=a_1 b_1+\cdots +a_t b_t$.  The following lemma
will be used later on.

\begin{lemma} \label{lem-pre} If $0\ne\b\in \z_n^t$, then
$\tsum{\a\in\z_n^t} q^{\a\cdot \b}=0.$
\end{lemma}
\begin{proof} Obviously, if $n \nmid j,$ then $\stackrel{n-1}{\tsum{i=0}}
q^{ij}=0$ since $q$ is an $n$-th root of 1. In general, if
$\b=(b_1, b_2, \cdots, b_t)\ne 0$, then there exists $i$ such that
$b_i\ne 0$. We have
\begin{eqnarray*}
\tsum{\a\in\z_n^t} q^{\a\cdot \b} =(\sum_{a_1 \in \z_n}
q^{a_1b_1})\cdots(\sum_{a_i \in \z_n} q^{a_ib_i})\cdots (\sum_{a_t
\in \z_n} q^{a_tb_t})=0.
\end{eqnarray*}
The lemma follows.
\end{proof}

Given a positive integer $m$ and a variable $x$, we denote
$$(m)_x=1+x+x^2+\cdots+x^{m-1},$$
and
$$(m)!_x=(m)_x (m-1)_x\cdots
(1)_x.$$
We set $(0)!_x=1$, $\qqbin{m }_x=\frac{(m)!_x}{(m-u)!_x (u)!_x}$ if
$u\leq m$ are positive integers. The notations $[m]_x, [0]!_x, [m]!_x,
$ and $\bin{m}{s}_x$ are as stated before.

A quiver $Q$ is an oriented graph given by two sets $Q_{0}$ and
$Q_{1}$ of vertices and arrows and two maps $s, t : Q_{1} \to
Q_{0}$ providing each arrow with its source and terminal vertex. A
path $\gamma$ is a finite sequence of concatenated arrows $\gamma=
a_{n}\cdots a_{1}$, which means that $t(a_{i}) = s(a_{i+1})$ for
$i = 1,\cdots,n-1$. We set $s(\gamma) = s(a_{1})$ and $t(\gamma) =
t(a_{n})$. It is noted that a vertex $u$ coincides with its source
and terminal vertex. The length of a path means the length of its
arrow sequence; vertices mean zero-length paths.

The path algebra of $Q$  is a $k$-vector space $kQ$ with the basis
$\{p |\ p$ is a path of $Q$\}, where the multiplication of two
paths $p$ and $q$ is defined by
$$p\cdot q=\left\{%
\begin{array}{ll}
pq, & \hbox{ if } s(p)=t(q),\\
    \\
0, &\hbox{otherwise.}\\
\end{array}%
\right.
 $$
Let $\ol{Q}$ be the double quiver of $Q.$  Namely, $\ol{Q}$ is
obtained by adding a reverse arrow $a^*:j \rightarrow i$ for each
arrow $a: i \rightarrow j$ in $Q.$ For each vertex $i,$ let $e_i$
be the associated trivial path.

We always assume that $\ell=n$ if $n$ is odd and $\ell=n/2$ if $n$
is even.

\section{Deformations of double quiver algebras}

For each pair $(\cc, \z_n)$, we associate to it a quiver
$\Q=\Q(\cc,\z_n)$ as follows.  The set of vertices has a
one-to-one correspondence with the group $\z_n^t.$ The set of
arrows is
$$\{ a(x,i): x-c^i \rightarrow x \ | \ x \in
Z_n^t, \ i \in \t \}.$$ It is noticed that $e_u \, a(u,
i)=a(u,i)\, e_{u-c^i}=a(u, i),\quad a(u, i)^*\, e_u=e_{u-c^i}\,
a(u, i)^*=a(u, i)^* $ in the double quiver algebra $k\ol{\Q}$.

The following lemma is well known.
\begin{lemma} {\rm (\cite{c2}, Proposition 2.3)} \label{lem2-1}
The path algebra $k\Q$ is a Hopf algebra where the
comultiplication, counit, and antipode are given by the following
\begin{eqnarray*}
\D(e_x) &=& \sum_{u+v=x} e_u \otimes e_v, \\
\D(a(x,i)) &=& \sum_{u+v=x}q^{u_i}e_u \otimes
a(v,i)+\sum_{u+v=x}a(u,i) \otimes e_v, \\
\e(e_0) &=& 1, \ \e(e_x)=0 \textrm{ for all} \ x \ne 0, \
\e(a(x,i))=0,\\
S(a(x, i)) &=& -q^{x_i-2} a(-x+c^i, i),\ S(e_x)= e_{-x},
\textrm{ for all } \ x, i.
\end{eqnarray*}
\end{lemma}

We consider the quotient algebra $\Pi^\cc=k\ol{\Q}/I,$
where $I$ is the ideal generated by the following relations
\begin{eqnarray}
 &&\label{eqn2-1} a(x,i)a(x,i)^*-a(x+c^i,i)^*a(x+c^i,i)-
    \frac{q^{x_i}-q^{-x_i}}{q-q^{-1}}e_x, \\
  &&\label{eqn2-2} a(x,j)^*a(x,i)-a(x-c^j,i)a(x-c^i,j)^*
\end{eqnarray}
for all $x\in\z_n^t,  i, j \in \t$ with  $i\ne j$.

The following fact is important to our aim.

\begin{lemma}
The algebra $\Pi^\cc$ is a Hopf algebra with comultiplication,
counit, and antipode given by
\begin{eqnarray*}
\D(e_x) &=& \sum_{u+v=x} e_u \otimes e_v, \\
\D(a(x,i)) &=& \sum_{u+v=x}q^{u_i}e_u \otimes
a(v,i)+\sum_{u+v=x}a(u,i) \otimes e_v, \\
\D(a(x,i)^*) &=& \sum_{u+v=x} e_u \otimes
a(v,i)^*+\sum_{u+v=x}q^{-v_i}a(u,i)^* \otimes e_v, \\
\e(e_0) &=& 1, \ \e(e_x)=0, \textrm {  for all } \ x \ne 0,\\
 \e(a(x,i))&=& 0, \ \ \e(a(x,i)^*)=0, \\
S(a(x, i)) &=& -q^{x_i-2} a(-x+c^i, i), S(a(x, i)^*)= -q^{-x_i+2}
a(-x+c^i, i)^*,\\
\quad S(e_x)&=&e_{-x}, \ \textrm{ for all } x\in\z_n^t, \  i\in\t.
\end{eqnarray*}

\end{lemma}

\begin{proof}
Let $\Q^{op}$  be the reverse quiver of $\Q$. Then we have a path
algebra $k\Q^{op}$. By duality and Lemma \ref{lem2-1}, one can
prove that $k\Q^{op}$ is  a Hopf algebra equipped with
\begin{eqnarray*}
\D(a(x,i)^*) &=& \sum_{u+v=x} e_u \otimes
a(v,i)^*+\sum_{u+v=x}q^{-v_i}a(u,i)^* \otimes e_v, \\
\e(e_0) &=& 1, \ \e(e_x)=0, \textrm {  for all } \ x \ne 0;\\
\e(a(x,i)^*)&=& 0, \\
 S(a(x, i)^*)&=& -q^{-x_i+2} a(-x+c^i, i)^*, \quad
\textrm{ for all  } \ x, i.\\
S(e_x)&=& e_{-x}, \textrm{ for  }  x\in \z_n^t.
\end{eqnarray*}
To see that $\Pi^\cc$ is  a Hopf algebra, it is sufficient to show
that the given maps $\Delta$ and $\e$ keep the relations (\ref{eqn2-1})
and (\ref{eqn2-2}). Indeed, for the relation (\ref{eqn2-1}), we
have
\begin{eqnarray*}
&&\D(a(x,i))\D(a(x,i)^*) -\D(a(x+c^i,i)^*)\D(a(x+c^i,i))\\
&=&\tsum{u+v=x} q^{u_i} e_u\otm a(v, i)a(v, i)^* +\tsum{u+v=x}
q^{-v_i} a(u, i) a(u, i)^* \otm e_v\\
&\quad &-\tsum{u+v=x+c^i} q^{u_i} e_u\otm a(v, i)^*a(v, i)-
\tsum{u+v=x+c^i} q^{-v_i} a(u, i)^* a(u, i)\otm e_v\\
&=&\tsum{u+v=x} q^{u_i} e_u\otm \left[ a(v, i)a(v, i)^*
-a(v+c^i, i)^*a(v+c^i, i)\right ]\\
&\quad & +\tsum{u+v=x} q^{-v_i}
\left[ a(u, i)a(u, i)^*-a(u+c^i, i)^*a(u+c^i, i)\right ]\otm e_v\\
&=&\tsum{u+v=x}q^{u_i} e_u\otm \frac{q^{v_i}-q^{-v_i}}{q-q^{-1}}
e_v+
\tsum{u+v=x} q^{-v_i} \frac{q^{u_i}-q^{-u_i}}{q-q^{-1}} e_u\otm e_v\\
&=&\left (\frac{q^{x_i}-q^{-x_i}}{q-q^{-1}}\right )\tsum{u+v=x} e_u\otm e_v.
\end{eqnarray*}
Hence $\D$ keeps the relation (\ref{eqn2-1}). For the relation
(\ref{eqn2-2}), it is easy to see that
\begin{eqnarray*}
&&\D(a(x,j)^*)\D(a(x,i)) -\D(a(x-c^j,i))\D(a(x-c^i,j)^*)\\
&=&\sum_{u+v=x} q^{u_i} e_u\otm \left[ a(v, j)^* a(v, i)-
a(v-c^j, i) a(v-c^i, j)^*\right ]\\
&\quad & + \sum_{u+v=x} q^{-v_i} \left [ a(u,j)^* a(u, i)-a(u-c^j,
i) a(u-c^i, j)^*\right ] \otm e_v
\\
&=&0
\end{eqnarray*}
The arguments for the counit $\e$ are similar. It remains to show
that $S$ is an antipode. For example,
\begin{eqnarray*}
&\quad& S(a(x,i)^*)S(a(x,i))-S(a(x+c^i,i))S(a(x+c^i,i)^*)\\
&=&a(-x+c^i, i)^* a(-x+c^i, i)-a(-x, i)a(-x, i)^*\\
&=&-\,
\frac{q^{-x_i}-q^{x_i}}{q-q^{-1}}e_{-x}=\frac{q^{x_i}-q^{-x_i}}{q-q^{-1}}S(e_x).
\end{eqnarray*}
For the relation (\ref{eqn2-2}) the argument for $S$ is similar
(we note that $a_{ij}=a_{ji}$). Hence $S$ is the antipode of
$\Pi^\cc$.  The lemma follows.
\end{proof}

We denote the paths $a(x, i_1)a(x-c^{i_1}, i_2) \cdots
a(x-c^{i_1}-\cdots-c^{i_{s-1}}, i_s)$ of   $\oq$  by $a(x, i_1i_2
\cdots i_s).$ For brevity, $a(x, i)a(x-c^i, i) \cdots
a(x-(s-1)c^i, i)$ is denoted by $a(x,i^s).$ Similarly, we use the
notations $a(x, i_1i_2 \cdots i_s)^*$ and $a(x, i^s)^*.$

\begin{lemma} \label{lem2-2}  We have the following formulae:
\begin{eqnarray*}
\D(a(x, i^m))&=&\sum_{u+v=x\atop s+t=m}\qbin{m}{s}_{q^{-2}}
q^{tu_i}
a(u, i^s)\otm a(v, i^t);\\
\D(a(x, i^m)^*)&=&\sum_{u+v=x\atop s+t=m}\qbin{m}{s}_{q^2} q^{-s
v_i} a(u, i^s)^*\otm a(v, i^t)^*.
\end{eqnarray*}
In particular, we have
\begin{eqnarray*}
\D(a(x, i^\ell))&=&\sum_{u+v=x} q^{\ell u_i} e_u\otm a(v,
i^\ell)+\sum_{u+v=x}
a(u, i^\ell)\otm e_v;\\
\D(a(x, i^\ell)^*)&=&\sum_{u+v=x} q^{-\ell v_i} a(u, i^\ell)^*\otm
e_v+\sum_{u+v=x} e_u\otm a(v, i^\ell)^*.
\end{eqnarray*}
\end{lemma}
\begin{proof} It is easy to see that the formula
\begin{eqnarray*}
\qqbin{m+1}_x=\qqbin{m}_x+x^{m-u+1}\left(\begin{array}{c}
    m\\
    u-1\\
    \end{array}\right )_x
\end{eqnarray*}
holds. By this formula, the proof is completed by induction
on $m$.
\end{proof}

We set $\kappa=1-a_{ij}$ and for all $x\in\z_n^t$,
$i,\; j\in\t$ with $i\ne j$,
  \begin{eqnarray*}
 \omega_{ij}(x)&=&
\sum_{t=0}^{\kappa}(-1)^t\bin{\kappa}{t}_q a(x,i^{\kappa-t}\,j\,i^t);\\
\omega_{ij}(x)^*&=& \sum_{t=0}^{\kappa}(-1)^t\bin{\kappa}{t}_q
a(x,i^{\kappa-t}\,j\,i^t)^*.
\end{eqnarray*}
%%%%%%%%%%%
%%%%%%%%%%%%
%%%%%%%%%%%%%%%
We have
\begin{lemma}\label{lem2-3}  For all $x\in\z_n^t$,
$i,\; j\in\t$ with $i\ne j$,
\begin{eqnarray*}
\D(\omega_{ij}(x))&=&\sum_{u+v=x} q^{\kappa\, u_i+u_j} e_u\otm
\omega_{ij}(v)+
\sum_{u+v=x} \omega_{ij}(u)\otm e_v\\
\D(\omega_{ij}(x)^*)&=&\sum_{u+v=x} e_u\otm \omega_{ij}(v)^*+
\sum_{u+v=x} q^{-\kappa\, v_i-v_j}\omega_{ij}(u)^*\otm e_v
\end{eqnarray*}
\end{lemma}
\begin{proof}
The proof is considerably  straightforward.
We compute the formula $\Delta(\omega_{ij}(x))$ where $a_{ij}=-1$.
In this case,
$$\omega_{ij}(x)=a(x, i^2j)-(q+q^{-1})a(x,iji)+a(x,j\,i^2).$$
We have
\begin{eqnarray*}
&\quad&\Delta(a(x, i^2j))=\\
&\quad& \sum_{u+v=x} q^{2u_i+u_j} e_u\otm a(v, i^2 j)
+\sum_{u+v=x} q^{2u_i} a(u, j)\otm a(v, i^2)\\
&\quad &+\sum_{u+v=x} (q+q^{-1}) q^{u_i+u_j} a(u, i)\otm a(v,
i\,j)\\
&\quad& \quad+\sum_{u+v=x} (1+q^{-2}) q^{u_i} a(u, i\,j)\otm a(v, i)\\
&\quad &\quad +\sum_{u+v=x} q^{u_j+2} a(u, i^2)\otm a(v, j)+\sum_{u+v=x}
a(u, i^2 j)\otm e_v.
\end{eqnarray*}
%%%%%%%%%
%%%%%%%%%%%
Similarly,
\begin{eqnarray*}
&\quad&\Delta(a(x, ji^2))=\\
&\quad& \sum_{u+v=x} q^{2u_i+u_j} e_u\otm a(v, j\;i^2)
+\sum_{u+v=x} q^{2u_i+2} a(u, j)\otm a(v, i^2)\\
&\quad &+\sum_{u+v=x} (1+q^{-2}) q^{u_i+u_j} a(u, i)\otm a(v,
j\,i)\\
&\quad& \quad +\sum_{u+v=x} (q+q^{-1}) q^{u_i} a(u, j\,i)\otm a(v, i)\\
&\quad &+\sum_{u+v=x} q^{u_j} a(u, i^2)\otm a(v, j)+\sum_{u+v=x}
a(u, j\;i^2)\otm e_v,
\end{eqnarray*}
and
\begin{eqnarray*}
&\quad&\Delta(a(x, i\,j\,i))=\\
&\quad& \sum_{u+v=x} q^{2u_i+u_j} e_u\otm a(v, i\,j\,i)
+\sum_{u+v=x} q^{u_i+u_j-1} a(u, i)\otm a(v, j\;i)\\
&\quad &+\sum_{u+v=x} q^{2u_i+1} a(u, j)\otm a(v, i^2)
+\sum_{u+v=x} q^{u_i-1} a(u, i\,j)\otm a(v, i)\\
&\quad &+\sum_{u+v=x} q^{u_i+u_j} a(u, i)\otm a(v, i\;j)+
\sum_{u+v=x} q^{u_j+1} a(u, i^2)\otm a(v, j)\\
&\quad &+ \sum_{u+v=x} q^{u_i} a(u, j\;i)\otm a(v, i)+
\sum_{u+v=x} a(u, i\;j\;i)\otm e_v.
\end{eqnarray*}
Therefore,
$$\D(\omega_{ij}(x))=\sum_{u+v=x} q^{2\, u_i+u_j} e_u\otm \omega_{ij}(v)+
\sum_{u+v=x} \omega_{ij}(u)\otm e_v.$$ The cases when $a_{ij}=0$
and $\omega_{ij}(x)^*$ are similar.
\end{proof}

Let $J$ be the ideal of $\Pi^\cc$ generated by
\begin{eqnarray}\label{eqn2-3}
 &&a(x,i^\ell), \ a(x,i^\ell)^*, \  \omega_{ij}(x),\ \omega_{ij}(x)^*
\end{eqnarray}
for all $x\in\z_n^t$ and
$i,\; j\in\t$ with $i\ne j$. By the definition of the
antipode, Lemmas \ref{lem2-2} and \ref{lem2-3}, it is easy to see
that $J$ is a Hopf ideal. We denote by $\uu^\cc$ the quotient
algebra $\Pi^\cc/J$, which can be presented
by generators and
relations as follows. As an algebra, $\uu^\cc$ is generated by
 $\{e_x, a(x, i), a(x, i)^*| x\in\z_n^t, i\in\t\}$
with the following relations: for all $x\in\z_n^t$ and
$i,\; j\in\t$ with $i\ne j$,
\begin{eqnarray*}
&& e_xe_y=\d_{xy}e_x, e_x \, a(y, i)=a(y,i)\, e_{x-c^i}=\d_{x, y} a(x, i),\\
&&a(x, i)^*\, e_y=e_{y-c^i}\, a(x, i)^*=\d_{x, y} a(x, i)^* ,\\
&&a(x,i)a(x,i)^*-a(x+c^i,i)^*a(x+c^i,i)=
    \frac{q^{x_i}-q^{-x_i}}{q-q^{-1}}e_x, \\
&&a(x,j)^*a(x,i)=a(x-c^j,i)a(x-c^i,j)^*,\\
&&a(x,i^\ell)=0, \ a(x,i^\ell)^*=0,\\
&&\omega_{ij}(x)=0, \ \omega_{ij}(x)^*=0.
\end{eqnarray*}
We yield the result as follows.
\begin{theorem}
$\uu^\cc$ is a Hopf algebra, of which the comultiplication, counit,
and antipode are defined by
\begin{eqnarray*}
\D(e_x) &=& \sum_{u+v=x} e_u \otimes e_v, \\
\D(a(x,i)) &=& \sum_{u+v=x}q^{u_i}e_u \otimes
a(v,i)+\sum_{u+v=x}a(u,i) \otimes e_v, \\
\D(a(x,i)^*) &=& \sum_{u+v=x} e_u \otimes
a(v,i)^*+\sum_{u+v=x}q^{-v_i}a(u,i)^* \otimes e_v, \\
\e(e_0) &=& 1, \ \e(e_x)=0, \textrm {  for all } \ x \ne 0,\\
 \e(a(x,i))&=& 0, \ \e(a(x,i)^*)=0,\
S(a(x, i)) = -q^{x_i-2} a(-x+c^i, i), \\
S(a(x, i)^*) &=& -q^{-x_i+2}
a(-x+c^i, i)^*, \textrm{ for all  } \ x, i,\\
S(e_x) &=& e_{-x},\  \textrm{ for all } x.
\end{eqnarray*}
\end{theorem}

\section{Isomorphisms between $\u$ and $\uu^\cc$}

In this section, we keep all notations as before.

For $x=(x_1, \cdots, x_t)$, let $K_x=K_1^{x_1}\cdots K_t^{x_t}$
and $\epsilon_x=\frac{1}{n^t} \sum\limits_{y \in \z_n^t}q^{-x\cdot
y} K_y.$  By Lemma \ref{lem-pre},
%%%%%%%
it is easy to  see that  for any $x, \ y\in \z_n^t$,
%%%%%%%%
\begin{eqnarray*}
 K _y \epsilon_x &=&q^{x\cdot y}
\frac{1}{n^t}\tsum{\b\in \z_n^t} q^{-x\cdot(\b+y)} K_{\b+y}
=q^{x\cdot y} \epsilon_x,\\
\epsilon_x \epsilon_y
&=&\frac{1}{n^t}\tsum{\gamma\in \z_n^t} q^{-x\cdot\gamma} K_\gamma
\epsilon_y= \frac{1}{n^t}\tsum{\gamma\in \z_n^t}
q^{-(x-y)\cdot\gamma} e_y=\delta_{x, y} \epsilon_x,\\
\sum\limits_{x\in \z_n^t} \epsilon_x &=&\sum\limits_{x\in \z_n^t}
\frac{1}{n^t}\tsum{\b\in \z_n^t} q^{-x\cdot\b}
K_{\b}=1+\frac{1}{n^t}\tsum{\b\ne 0\atop{\b\in \z_n^t}}
(\sum\limits_{x\in \z_n^t} q^{-x\cdot\b}) K_{\b}=1,\\
%$$
%It is easy to see that for all $x\in\z_n^t$,
E_i \epsilon_x &=& \epsilon_{x+c^i} E_i, \  F_i \epsilon_x=\epsilon_{x-c^i} F_i,\\
\D(\epsilon_x) &=& \sum_{u+v=x} \epsilon_u\otm\epsilon_v, \
S(\epsilon_x)=\epsilon_{-x}.
\end{eqnarray*}

We have constructed the Hopf algebra $\uu^\cc$. The relationship between
$\uu^\cc$ and $\u$ is given as follows.
\begin{theorem} \label{thm2} The map
$\tau: \uu^\cc\to \u$ defined by
$$%\begin{eqnarray*}
\tau(e_x)=\epsilon_x,  \
\tau(a(x, i))= \epsilon_x\, E_i, \
\tau(a(x, i)^*)= F_i\,\epsilon_x
$$%\end{eqnarray*}
is a Hopf algebra isomorphism.
\end{theorem}

\begin{proof} We prove the theorem in several steps.
\vskip 0.3truecm

{\it Step 1}: The map $\tau$ is well-defined. We should first
verify that $\tau$ keeps the basic relations for path algebra.
That is, $e_xe_y=\d_{xy}e_x,$ $e_xa(y,i)=a(y, i)e_{x-c^i}=\d_{xy}a(y,i)$,
and
$a(x,i)^*e_y=e_{y-c^i} a(x, i)^*=\d_{xy}a(x,i)^*.$  For the first one,
\begin{eqnarray*}
\tau(e_x)\tau(e_y) &=& \epsilon_x\epsilon_y=\delta_{x, y}
\epsilon_x=\tau(\d_{xy}e_x) =\tau(e_x\,e_y).
\end{eqnarray*}
For the second one,
\begin{eqnarray*}
\tau(e_x)\tau(a(y, i))=\epsilon_x\epsilon_y E_i=\d_{x, y}
\epsilon_y E_i=\tau(\d_{xy}\, a(y, i)).
%=\epsilon_y E\epsilon_{x-c^i}=\tau(a(y, i))\tau(e_{x-c^i}).
\end{eqnarray*}
The rest are similar.

It remains to show that $\tau$ keeps the relations
(\ref{eqn2-1})-(\ref{eqn2-3}). For example, for the relation (\ref{eqn2-1}),
\begin{eqnarray*}
&&\tau(a(x,i))\tau(a(x,i)^*)-\tau(a(x+c^i,i)^*)\tau(a(x+c^i,i))-
    \frac{q^{x_i}-q^{-x_i}}{q-q^{-1}}\tau(e_x)\\
    &&=\epsilon_x E_i F_i\epsilon_x- (F_i\epsilon_{x+c^i})
    (\epsilon_{x+c^i} E_i)-
\frac{q^{x_i}-q^{-x_i}}{q-q^{-1}}\epsilon_x\\
    &&=(E_iF_i-F_i E_i)\epsilon_x-\frac{q^{x_i}-q^{-x_i}}{q-q^{-1}}\epsilon_x\\
    &&=\frac{K_i-K_i^{-1}}{q-q^{-1}}\epsilon_x-\frac{q^{x_i}-q^{-x_i}}{q-q^{-1}}\epsilon_x\\
&&=\frac{q^{x_i}-q^{-x_i}}{q-q^{-1}}\epsilon_x-
    \frac{q^{x_i}-q^{-x_i}}{q-q^{-1}}\epsilon_x=0.
\end{eqnarray*}
For the relation (\ref{eqn2-2}), it is similar. For the relation
(\ref{eqn2-3}),
$$\tau(a(x,i))\cdots \tau(a(x-(\ell-1)c^i, i)=
(\epsilon_x E_i)\cdots (\epsilon_{x-(\ell-1)c^i} E_i)=\epsilon_x
E_i^\ell=0,
$$
and
\begin{eqnarray*}
&\quad & \sum_{t=0}^{\kappa}(-1)^t\bin{\kappa}{t}_q
\tau(a(x,i))\cdots
\tau(a(x-(\kappa-t-1)c^i, i)\\
&\quad&\quad\times\tau(a(x-(\kappa-t)c^i, j))\\
&\quad& \quad \times \tau(a(x-(\kappa-t)c^i-c^j), i)
\cdots \tau(a(x-(\kappa-1)c^i-c^j, i))\\
&=&\sum_{t=0}^{\kappa}(-1)^t\bin{\kappa}{t}_q (\epsilon_x
E_i)\cdots (\epsilon_{x-(\kappa-t-1)c^i} E_i)\\
&\quad&\quad\times (\epsilon_{x-(\kappa-t-1)c^i-c^j} E_i)
(\epsilon_{x-(\kappa-t)c^i} E_j)\cdots ( \epsilon_{x-(\kappa-1)c^i-c^j} E_i)\\
&=&\epsilon_x\bigg(\sum_{t=0}^{\kappa}(-1)^t\bin{\kappa}{t}_q
E_i^{\kappa-t} E_j E_i^t\bigg)=0.
\end{eqnarray*}
The rest relations in (\ref{eqn2-3}) are similar.

\vskip 0.3truecm

{\it Step 2}:  Define an algebra map $\sigma: \u\to\uu^\cc$ by
$$
\sigma(K_c)= \sum_{x \in \z_n^t}q^{c\cdot x} e_x,\
\sigma(E_i)=\sum_{x \in \z_n^t} a(x,i),\ \sigma(F_i)= \sum_{x \in
\z_n^t} a(x,i)^*. $$
The aim is to show that $\sigma$ is the
inverse of $\tau$. The map $\sigma$ is also well defined. Indeed,
for the relation (\ref{eqn3-1}),
$$ \sigma(K_i)^n=\bigg(\s q^{x_i}e_x\bigg)^n=\sum_{x\in\z_n^t} q^{nx_i}e_x
=\sum_{x\in\z_n^t} e_x=1,$$ since $q$ is an $n$-th root of unity.
For the relation (\ref{eqn3-2}), we have
\begin{eqnarray*}
\sigma(K_i)\sigma(E_j)\sigma(K_i^{-1}) &=&\bigg (\s q^{x_i}e_x
\bigg)\bigg(\s a(x,j)\bigg)
\bigg(\s q^{-x^i}e_x\bigg )\\
&=& \s q^{x_i} e_x a(x,j)\; q^{-x_i+a_{ij}}\; e_{x-c^j} \\
&=& q^{a_{ij}}\s a(x,j)=q^{a_{ij}}\sigma(E_j).
\end{eqnarray*}
The another relation in (\ref{eqn3-2}) is similar. For the
relation (\ref{eqn3-3}),
\begin{eqnarray*}
&\quad& \sigma(E_i)\sigma(F_j)-\sigma(F_j)\sigma(E_i)\\
&=&\bigg(\s a(x,i)\bigg )\bigg (\s a(x,j)^*\bigg)
-\bigg(\s a(x,j)^*\bigg)\bigg(\s a(x,i)\bigg )\\
%&=&\s a(x,i)a(x,j)^*-\s a(x,j)^*a(x,i) \\
&=& \left\{%
\begin{array}{ll}
   \s \big [a(x,i)a(x,i)^*-a(x+c^i,i)^*a(x+c^i,i)\big], & \hbox{if $i=j$;} \\
   \\
   -\s \big[a(x,j)^*a(x,i)-a(x-c^j,i)a(x-c^i,j)^*\big], & \hbox{if $i \ne j$} \\
\end{array}%
\right. \\
&=& \left\{%
\begin{array}{ll}
   \s \frac{q^{x_i}-q^{-x_i}}{q-q^{-1}}e_x, & \hbox{if $i=j$;} \\
   \\
   0, & \hbox{if $i \ne j$} \\
\end{array}%
\right.\\
&=&\sigma\bigg(\delta_{ij}\frac{K_i-K_i^{-1}}{q-q^{-1}}\bigg).
\end{eqnarray*}
The relation (\ref{eqn3-4}) is due to the the first two relations
of (\ref{eqn2-3}). As for the quantum Serre relations, if
$a_{ij}=0,$ we have
\begin{eqnarray*}
&\quad& \sigma(E_i)\sigma(E_j)-\sigma(E_j)\sigma(E_i)\\
&=& \bigg(\s a(x,i)\bigg )\bigg(\s a(x,j)\bigg )-\bigg (\s
a(x,j)\bigg)
\bigg(\s a(x,i)\bigg )\\
&=& \s [a(x,ij)-a(x,ji)]=0;
\end{eqnarray*}
if $a_{ij}=-1,$ we have
\begin{eqnarray*}
& &\sigma(E_i)^2\sigma
(E_j)-(q+q^{-1})\sigma(E_i)\sigma(E_j)\sigma(E_i)
+\sigma(E_j)\sigma(E_i)^2\\
&=&\bigg(\s a(x,i)\bigg )^2\bigg (\s a(x,j)\bigg)-
(q+q^{-1})\bigg(\s a(x,i)\bigg)\\
&\quad&\times \bigg(\s a(x,j)\bigg )\bigg(\s a(x,i)\bigg)
+\bigg(\s a(x,j)\bigg )\bigg(\s a(x,i)\bigg)^2\\
&=& \s \omega_{ij}(x)=0.
\end{eqnarray*}
The arguments of the rest relations are similar.

By direct calculation, we have $\sigma \circ \tau=1$ and $\tau
\circ \sigma=1.$  Hence $\tau$ is an algebra isomorphism.

\vskip 0.3truecm

{\it Step 3}: $\tau$ is a Hopf algebra homomorphism. It is enough
to verify that $\tau$ is also a coalgebra map since bialgebra
homomorphisms are Hopf homomorphisms. It suffices to check it on
the generators, but this is considerably direct. We check only one
of them.
\begin{eqnarray*}
\D(\tau(a(i, x)))&=&\D(\epsilon_x E_i)=\D(\epsilon_x)\,\D(E_i)\\
&=&\bigg(\sum_{u+v=x}\epsilon_u\otm\epsilon_v\bigg )
\bigg(E_i\otm 1+K_i\otm E_i\bigg )\\
&=&\sum_{u+v=x} \tau(a(u, i))\otm \tau(e_v)+\sum_{u+v=x}
K_i\epsilon_u\otm \tau(a(v, i))\\
&=&\sum_{u+v=x} \tau(a(u, i))\otm \tau(e_v)+\sum_{u+v=x}
q^{u_i}\tau(e_u)\otm \tau(a(v, i))\\
&=&(\tau\otm\tau)\D(a(x, i)),
\end{eqnarray*}
The others are similar.

The proof is completed.
\end{proof}

For the Cartan matrix $\cc=(a_{ij})_{t\times t}$, there is an
associated Hopf algebra $\uu^+$ generated by $K_i, \ E_i$ for $1
\le i \le t$, subjecting to the relations
\begin{eqnarray*}
     &&K_i^n=1, \ K_iK_j=K_jK_i,\\
    &&K_iE_jK_i^{-1}=q^{a_{ij}}E_j.\\
&&\D(K_i)=K_i \otimes K_i, \\
&&\D(E_i)=E_i \otimes 1 + K_i
\otimes E_i.\\
&&\e(K_i)=1, \ \e(E_i)=0,\\
&&S(K_i)=K_i^{-1}, \ S(E_i)=-K_i^{-1} E_i.
\end{eqnarray*}
Similarly, there is an associated Hopf algebra $\uu^-$ generated
by $K_i, \ F_i$ for $1 \le i \le t$, subjecting to the relations
\begin{eqnarray*}
     &&K_i^n=1, \ K_iK_j=K_jK_i; \\
    &&K_iF_jK_i^{-1}=q^{-a_{ij}}F_j;\\
&&\D(K_i)=K_i \otimes K_i, \\
&& \D(F_i)=F_i \otimes K_i^{-1} + 1 \otimes F_i,\\
&&\e(K_i)=1, \ \ \e(F_i)=0,\\
&&S(K_i)=K_i^{-1},  \ S(F_i)=-F_i K_i.
\end{eqnarray*}
By \cite{x}, there is a skew Hopf pairing $\uu^+\times \uu^-\to k$
and  we have the Drinfeld  double ${\mathcal{D}}(\uu)$ (see
\cite{x}, Sect. 2). Let
$I=\langle K_i\otm 1-1\otm K_i\ |\ i\in\t\rangle$ be the ideal of
${\mathcal{D}}(\uu)$ generated
by $K_i\otm 1-1\otm K_i, \ i\in\t$.
It is easy to see that $I$ is a Hopf ideal of ${\mathcal{D}}(\uu)$
and $\Pi^\cc\cong {{\mathcal{D}}(\uu)}/{I}$
as Hopf algebras by Theorem \ref{thm2}. Furthermore, if $\t=\{1\}$,
the relation (\ref{eqn2-2}) automatically vanishes and $\Pi^\cc$ is just
a deformation  of preprojective algebra associated to the quiver $\Q$.

 We have found an algebraic realization of the
quantum group $\u$. This method is very intuitive. It is natural
to expect that the presentation via a double quiver will help to
study the representation theory, probably by consulting the theory
of deformed preprojective algebras. This will be considered in the
forthcoming papers.

\section{Some Remarks}
In the representation theory of finite dimensional algebras,
finite dimensional basic algebras can  always be constructed
via quivers with
admissible relations, according to Gabriel's Theorem. Recall that
a path relation is called admissible if the length of the paths
involved are at least two. We remark that our relation
(\ref{eqn2-1}) is not admissible. Actually, the quantum groups
$\u$ are not basic and hence there is no hope to present them via
quivers with admissible relations. However it will be of interest
to consider the Ext-quiver of the quantum groups and compare with
the double quivers obtained.

For ${\frak{sl}}_2$ the quiver obtained by  the authors is the
same as the quiver described in arXiv: math.RT/0410017, which
appeared after the submission of our paper. The authors thank the
referee for pointing out the reference.

\section*{Acknowledgment}

\quad The authors are grateful to the referee for his/her helpful
comments and suggestions.


\begin{thebibliography}{}

\bibitem{c1} \label{c1} C. Cibils, {\it A quiver quantum group},
 Comm. Math.\ Phys. 157 (1993) 459-477.

\bibitem{c2}\label{c2}
C. Cibils, {\it Half-quantum groups at roots of unity, path
algebras, and representation type}, Inter. Math. Res. Noticies 12
(1997) 541-553.

\bibitem {cbh}\label{cbh} W. Crawley-Boevey, M. P. Holland,
{\it Noncommutative deformations of Kleinian Singularities}, Duke
Math.  J.  92 (1998) 605-635.

\bibitem{egst} Erdmann, K., Green, E. L., Snashall, N., Taillefer,
R. {\it Representation theory of the Drinfeld doubles of a family
of Hopf algebras}, preprint, arXiv: math.RT/0410017.


\bibitem{g} \label{g}  J. A. Green,  {\it Hall algebras, hereditary
algebras and quantum groups}, Invent. Math. 120 (1995) 361-377.

\bibitem{ig} \label{ig} I. Gordon, {\it Quantised
function algebras at roots of unity and path algebras}, J. Algebra
220 (1999) 381-395.

\bibitem{l1}
G. Lusztig, {\it Modular representations and quantum groups},
Contemp. Math. 82 (1989) 59-77.

\bibitem{l2}\label{l2}
G. Lusztig, {\it Finite dimensional Hopf algebras arising from
quantized universal enveloping algebras}, Journal of the AMS 3
(1990) 257-296.

\bibitem{l3}
G. Lusztig, {\it Quantum groups at roots of 1}, Geom. Dedicata 35
(1990) 89-113.


\bibitem{r} \label{r}  C. M. Ringel, { \it Hall algebras and quantum
groups}, Invent. Math. 101 (1990) 583-592.

\bibitem{x} \label{x}  J. Xiao, {\it Drinfeld double and Ringel-Green
theory of Hall Algebras},  J. Algebra 190 (1997) 100-144.

\bibitem{yang} \label{yang} Shilin Yang, {\it
Quantum groups and deformations of preprojective algebras},
 J. Algebra, 279 (2004), 3-21.
\end{thebibliography}
\end{document}